\documentclass[a4paper,11pt]{article}
 \usepackage[plainpages=false]{hyperref}
 \usepackage{amsfonts,amsmath,amssymb,amsthm}
 \usepackage{latexsym,lscape,rawfonts,mathrsfs}

 \usepackage[dvips]{color}
 \usepackage{multicol}

 \usepackage[all]{xy}
 \usepackage{eufrak}
 \usepackage{makeidx}
 \usepackage{graphicx,psfrag}

 \usepackage{array,tabularx}

 \usepackage{setspace}

 \usepackage{appendix}

\usepackage{txfonts}

 \newcommand{\C}{\ensuremath{\mathbb{C}}}

 \newcommand{\ba}{\begin{align*}}
 \newcommand{\ea}{\end{align*}}

 \makeatletter
 \def\ExtendSymbol#1#2#3#4#5{\ext@arrow 0099{\arrowfill@#1#2#3}{#4}{#5}}
 
 \makeatother

 \makeatletter
 \def\ExtendSymbol#1#2#3#4#5{\ext@arrow 0099{\arrowfill@#1#2#3}{#4}{#5}}
 
 \makeatother

 \definecolor{hao}{rgb}{1,0.5,0}
 \definecolor{miao}{cmyk}{0.5,0,0.2,0.2}
 \definecolor{qiao}{gray}{0.96}

 \title{Two Comments on ``Laplace comparison on K\"ahler Ricci flow and convergence (\href{https://arxiv.org/pdf/2509.14820}{arXiv:2509.14820v1. })"}
 \author{Xiuxiong Chen, Bing Wang}
 \date{Oct 30, 2025}

\begin{document}
 \maketitle
 
{\bf (1).} In the concluding acknowledgments of the recent preprint~\cite{TZZZZ}, one of the authors (Prof. Q.~S.~Zhang) thanks us for ``helpful comments''. 
 This is confusing, as we had no opportunity to read or comment on the preprint prior to its posting. We note only that, more than six years ago, we exchanged emails with the same author addressing questions about our earlier paper.
 
{\bf (2).} On pages 40-41 of~\cite{TZZZZ}, the authors challenge the definition of the \emph{polarized canonical radius} (pcr) in~\cite{CW4}, claiming--based on a “direct computation”--that the pcr cannot be bounded from below. Their claim, in essence, is the following (lines 8-21, page 41 of~\cite{TZZZZ}):
 
 \textit{
 Let $(M, g, J, L,h)$ be a polarized K\"ahler manifold.  Then we proceed by rescaling. 
 Let $\tilde{g}=4kj_0^2 g$ and $\{\tilde{S}_i^{(j)}\}$ be orthonormal basis of $H^0(M, L^j)$ under the metrics $\tilde{\omega}$ and $\tilde{h}$.   Then
 \begin{align*}
   \inf_M \sum_{i=0}^{N_j} ||\tilde{S}_i^{(j)}||_{\tilde{h}(t)}^2(x,t) \leq \frac{1}{|M|_{\tilde{g}}} \int_M \sum_{i=0}^{N_j} ||\tilde{S}_i^{(j)}||_{\tilde{h}(t)}^2(x,t) d\tilde{g} \leq \frac{CN_j}{|M|_{\tilde{g}}}=\frac{CN_j}{(4kj_0^2)^{\frac{n}{2}}|M|_{g}} \to 0, 
 \end{align*}
 as $k \to \infty$.  Thus, $\inf_M \tilde{b}^{(j)}(x,t) \to -\infty$. 
 }
 
The \textbf{error} in their argument occurs on line 14 of page 41 in reference~\cite{TZZZZ}, where the authors assume that the line bundle remains unchanged when working on the manifold $(M, \tilde{g})$. However, to satisfy the polarization condition, the line bundle must be appropriately adjusted to $\tilde{L}$. Specifically, if $\tilde{g}=4kj_0^2 g$, then $\tilde{L}=L^{4kj_0^2}$. Accordingly, $\{\tilde{S}_i^{(j)}\}$ should be taken as an orthonormal basis of
$H^0(M, \tilde{L}^j)=H^0(M, L^{4kj_0^2 j})$
 with respect to $d\tilde{g}=(4kj_0^2)^{m} dg$ and $\tilde{h}^j=h^{4kj_0^2j}$. 
 Define 
 \begin{align*}
     \tilde N_j := \mathrm{dim}\, H^0(M, L^{4kj_0^2 j}) -1. 
 \end{align*}
 Running their argument with this correction yields
 \begin{align*}
   \inf_M \sum_{i=0}^{\tilde N_j} ||\tilde{S}_i^{(j)}||_{\tilde{h}^j(t)}^2(x,t) \leq \frac{1}{|M|_{\tilde{g}}} \int_M \sum_{i=0}^{\tilde N_j} ||\tilde{S}_i^{(j)}||_{\tilde{h}^j(t)}^2(x,t) d\tilde{g} \leq \frac{C \tilde N_j}{(4kj_0^2)^{\frac{n}{2}}|M|_{g}}.  
 \end{align*}
 However,  as $k \to \infty$, by Riemann-Roch theorem ($m=\frac{n}{2}=\dim_{\C} M$), we have
 \begin{align*}
        \tilde{N}_j \sim  \frac{c_1^{m}(L)}{m!} (4k j_0^2j)^m + O(k^{m-1})=|M|_g (4k j_0^2j)^{\frac{n}{2}} + O(k^{\frac{n}{2}-1}). 
 \end{align*}
 It follows that
 \begin{align*}
    \lim_{k \to \infty} \frac{C \tilde N_j}{(4kj_0^2)^{\frac{n}{2}}|M|_{g}} = C j^{\frac{n}{2}}>0. 
 \end{align*}
 Therefore, one cannot conclude that 
 \begin{align*}
   \inf_M \sum_{i=0}^{\tilde N_j} ||\tilde{S}_i^{(j)}||_{\tilde{h}^j(t)}^2(x,t) \to 0, \quad \textrm{as} \; k \to \infty, 
 \end{align*}
and consequently one cannot conclude that
  \begin{align*}
  \inf_M \tilde{b}^{(j)}(x,t) = \inf_M \log \sum_{i=0}^{\tilde N_j} ||\tilde{S}_i^{(j)}||_{\tilde{h}^j(t)}^2(x,t) \to -\infty, \quad \textrm{as} \; k \to \infty,  
 \end{align*}
as claimed in~\cite{TZZZZ}, line 22 on page 41.\\

Actually, one can show that the conclusion $\inf_M \tilde{b}^{(j)}(x,t) \to -\infty$ is impossible. 
Since $\{\tilde{S}_i^{(j)}\}$ is  an orthonormal basis of $H^0(M, \tilde{L}^j)=H^0(M, L^{4kj_0^2 j})$
with respect to $d\tilde{g}=(4kj_0^2)^{m} dg$ and $\tilde{h}^j=h^{4kj_0^2j}$,  it is clear that $\{(4kj_0^2)^{\frac{m}{2}}\tilde{S}_i^{(j)}\}$ is an orthonormal basis of $H^0(M, \tilde{L}^j)$ with respect to $dg$ and $\tilde{h}^j$. 
By the expansion formula of Bergman kernel (Yau-Tian-Zelditch-Lu-Catlin, cf.~\cite{Zelditch}), we have 
 \begin{align*}
      \sum_{i=0}^{\tilde N_j} ||(4kj_0^2)^{\frac{m}{2}} \tilde{S}_i^{(j)}||_{\tilde{h}^j(t)}^2(x,t)=  \frac{(4kj_0^2 j)^{m}}{\pi^m} + O(k^{m-1})
 \end{align*}
 which is the same as
  \begin{align*}
      \sum_{i=0}^{\tilde N_j} ||\tilde{S}_i^{(j)}||_{\tilde{h}^j(t)}^2(x,t)=\frac{j^m}{\pi^m} + O(k^{-1}). 
 \end{align*}
 Therefore, for each fixed $x \in M$, we have
  \begin{align*}
    \tilde{b}^{(j)}(x,t)=\log   \sum_{i=0}^{\tilde N_j} ||\tilde{S}_i^{(j)}||_{\tilde{h}^j(t)}^2(x,t)=\log \frac{j^m}{\pi^m} + O(k^{-1}) \to \log \frac{j^m}{\pi^m}  \neq -\infty, \quad \textrm{as} \; k \to \infty. 
 \end{align*}
 Since $M$ is compact, the above asymptotic behavior excludes the possibility: $\inf_M  \tilde{b}^{(j)}(x,t) \to -\infty$ as $k \to \infty$. \\

 \textbf{Updates:}   The authors of~\cite{TZZZZ} recently update their paper and move the appendix to a new note~\cite{TZZZZ2}, where they basically repeat their so called ``counter-example".
 Correspondingly, we put a further response in our new note~\cite{CWd}.

 \end{document}